\documentclass[leqno,12pt]{article}
\usepackage{amssymb,amsfonts}
\usepackage{amsmath,latexsym}

\usepackage{amstext,epic,eepic,epsf,pslatex}
\usepackage{graphicx}

\newcommand{\R}{{\mathbb R}}

\def\XXint#1#2#3{{\setbox0=\hbox{$#1{#2#3}{\int}$ }
\vcenter{\hbox{$#2#3$ }}\kern-.6\wd0}}

\newtheorem{thm}{Theorem}[section]

\begin{document}

\title{Multi-dimensional scalar conservation laws with unbounded integrable initial data}

\author{Denis Serre \\ \'Ecole Normale Sup\'erieure de Lyon\thanks{U.M.P.A., UMR CNRS--ENSL \# 5669. 46 all\'ee d'Italie, 69364 Lyon cedex 07. France. {\tt denis.serre@ens-lyon.fr}}}

\date{}

\maketitle

\begin{abstract}
We discuss the minimal integrability needed for the initial data, in order that the Cauchy problem for a multi-dimensional conservation law admit an entropy solution. In particular we allow unbounded initial data. We investigate also the decay of the solution as time increases, in relation with the nonlinearity.

The main ingredient is our recent theory of divergence-free positive symmetric tensor. We apply in particular the so-called {\em compensated integrability} to a tensor which generalizes the one that L. Tartar  used in one space dimension. It allows us to establish a Strichartz-like inequality, in a quasilinear context.

This program is carried out in details for a multi-dimensional version of the Burgers equation. 
\end{abstract}

\paragraph{Notations.} When $1\le p\le\infty$, the natural norm in $L^p(\R^n)$ is denoted $\|\cdot\|_p$, and the conjugate exponent of $p$ is $p'$. The total space-time dimension is $d=1+n$ and the coordinates are $x=(t,y)$. In the space of test functions, ${\cal D}^+(\R^{1+n})$ is the cone of functions which take non-negative values.  The partial derivative with respect to the coordinate $y_j$ is $\partial_j$, while the time derivative is $\partial_t$. Various finite positive constants that depend only the dimension, but not upon the solutions of our PDE, are denoted $c_d$~; they usually differ from one inequality to another one. 

\section{Introduction}

Let us consider a scalar conservation law in $1+n$ dimensions
\begin{equation}
\label{eq:mdcl}
\partial_tu+\sum_{i=1}^n\partial_if_i(u)=0,\qquad t>0,\,y\in\R^n.
\end{equation}
We complement this equation with an initial data
$$u(0,y)=u_0(y),\qquad y\in\R^n.$$
The flux $f(s)=(f_1(s),\ldots,f_n(s))$ is  a smooth vector-valued function of $s\in\R$. We recall the terminology that an entropy-entropy flux pair is a couple $(\eta,q)$ where $s\mapsto\eta(s)$ is a numerical function, $s\mapsto q(s)$ a vector-valued function, such that $q'(s)\equiv\eta'(s)f'(s)$. The Kruzhkov's entropies and their fluxes form a one-parameter family:
$$\eta_a(s)=|s-a|,\qquad q_a(s)={\rm sgn}(u-a)\,(f(u)-f(a)).$$
Together with the affine functions, they span the cone of convex functions. 

\bigskip

We recall that an {\em entropy solution} is a measurable function $u\in L^1_{\rm loc}([0,+\infty)\times\R^n)$ such that $f(u)\in L^1_{\rm loc}([0,+\infty)\times\R^n)$, which satisfies the Cauchy problem in the distributional sense,
\begin{equation}
\label{eq:distr}
\int_0^\infty dt\int_{\R^n}(u\partial_t\phi+f(u)\cdot\nabla_y\phi)\,dy+\int_{\R^n}u_0(y)\phi(0,y)\,dy=0,\qquad\forall\phi\in{\cal D}(\R^{1+n}),
\end{equation}
together with the entropy inequalities
\begin{eqnarray}
\nonumber
\int_0^\infty dt\int_{\R^n}(\eta_a(u)\partial_t\phi+q_a(u)\cdot\nabla_y\phi)\,dy & & \\
\label{eq:entrina}
+\int_{\R^n}\eta_a(u_0(y))\phi(0,y)\,dy & \ge0, & \qquad\forall\,\phi\in{\cal D}^+(\R^{1+n}),\,\forall\,a\in\R.
\end{eqnarray}
When it enjoys higher integrability, an entropy solution is expected to satisfy additional entropy inequalities of the form
\begin{equation}
\label{eq:entrin}
\int_0^\infty dt\int_{\R^n}(\eta(u)\partial_t\phi+q(u)\cdot\nabla_y\phi)\,dy+\int_{\R^n}\eta(u_0(y))\phi(0,y)\,dy\ge0,\qquad\forall\phi\in{\cal D}^+(\R^{1+n}),
\end{equation}
for more general convex entropies $\eta$. In particular, one is interested in inequality (\ref{eq:entrin}) for the entropy-entropy flux pair
$$\bar\eta(s)=\frac12\,s^2,\qquad\bar q(s)=\int_0^szf'(z)\,dz.$$

\bigskip

The theory of this Cauchy problem dates back to 1970, when S. Kruzhkov \cite{Kru} proved that if $u_0\in L^\infty(\R^n)$, then there exists one and only one entropy solution in the class
$$L^\infty(\R_+\times\R^n)\cap C(\R_+;L^1_{\rm loc}(\R^n)).$$
The operator $S_t:u_0\mapsto u(t,\cdot)$, which maps $L^\infty(\R^n)$ into itself, enjoys several additional properties. On the one hand, a comparison principle says that if $u_0\le v_0$, then $S_tu_0\le S_tv_0$. For instance, the solution $u$ associated with the data $u_0$ is majorized by the solution $\bar u$ associated with the data $(u_0)_+$, the positive part of $u_0$. On another hand, if $v_0-u_0$ is integrable over $\R^n$, then $S_tv_0-S_tu_0$ is integrable too, and
\begin{equation}
\label{eq:contr}
\int_{\R^n}|S_tv_0-S_tu_0|(y)\,dy\le\int_{\R^n}|v_0-u_0|(y)\,dy.
\end{equation}
Finally, if $u_0$ belongs to some $L^p(\R^n)$ space, then $S_tu_0$ has the same integrability, and the map $t\mapsto \|S_tu_0\|_p$ is non-increasing.
We warn the reader  that the contraction property (\ref{eq:contr}) occurs only for the $L^1$-norm, but not for other $L^p$-norms.

\bigskip

Because of (\ref{eq:contr}) and the density of $L^1\cap L^\infty(\R^n)$ in $L^1(\R^n)$, the family $(S_t)_{t\ge0}$ extends in a unique way as a continuous semi-group of contractions over $L^1(\R^n)$, still denoted $(S_t)_{t\ge0}$. When $u_0\in L^1(\R^n)$ is unbounded, we are thus tempted to declare that $u(t,y):=(S_tu_0)(y)$ is the {\em abstract solution} of the Cauchy problem for (\ref{eq:mdcl}) with initial data $u_0$. An alternate construction of $(S_t)_{t\ge0}$, based upon the Generation Theorem for nonlinear semigroups, is due to M. Crandall \cite{Cra}, who pointed out that  it is unclear  whether $u$ is an entropy solution, because the local integrability of the flux $f(u)$ is not guaranted. It is therefore an important question to identify the widest class of integrable data for which $u$ is actually an entropy solution of (\ref{eq:mdcl}). 

\bigskip

To achieve this goal, we develop a new strategy, based on the Compensated Integrability that we introduced in our previous papers \cite{Ser_DPT,Ser_CI}. It uses a map $a\mapsto M(a)\in{\bf Sym}_d$, whose lines are entropy-entropy flux pairs, where the entropies are precisely the functions ${\rm id}_\R,f_1,\ldots,f_n$ which appear in the conservation law. The map $M$ is a non-decreasing function of $a$. This tensor was already used when $n=1$ by L. Tartar \cite{Tar_HW} to prove the compactness of the semi-group, and by F. Golse \cite{Gol} (see also \cite{GoPe}) to prove some kind of regularity. An essential ingredient is the amount of non-linearity displayed by the flux $f$. We illustrate our strategy by
carrying out the details on the most typical nonlinear conservation law, a multi-d generalization of the Burgers equation.

\paragraph{Outline of the article.}
We begin with a detailed, definitive, analysis of the multi-d Burgers equation. The equation is described in the next section. Our main result is a well-posedness when the initial data is integrable. It is based on a dispersion estimate, which has the flavour of a Strichartz inequality, from which we derive a decay estimate of $L^p$-norms for $p\le\frac{d^2}{d-1}$\,. The proof is given in Sections \ref{s:disp} and \ref{s:dec}. We explain how the strategy extends to general fluxes $f$ in Section \ref{s:gen}.

\paragraph{Acknowledgements.} I am indebted to C. Dafermos, whose precious comments helped me to improve this article, in particular in giving full credit to previous contributors. I also thank L. Silvestre for correcting a miscalculation.

\section{The multi-d Burgers equation}\label{s:Bur}

For a conservation law of the general form (\ref{eq:mdcl}), it is harmless to assume $f(0)=0$. By chosing an appropriate inertial frame, which does not affect the norms $\|u(t)\|_p$, we may also assume $f'(0)=0$. Thus $f(s)=O(s^2)$ at the origin. Say that $f(s)\sim s^k\vec v_1$ as $s\rightarrow0$, where $\vec v_1$ is a non-zero constant vector. We may perform a linear change of the spacial coordinates such that $f_1(s)\sim \frac{s^k}k$\, and $f_j(s)=o(|s|^k)$ otherwise. Unless we meet a flat component, he process can be continued until we find a new coordinate system $(y_1,\ldots,y_n)$ in which
$$f_j(s)\sim\frac{s^{k_j}}{k_j}\,,\qquad 2\le k_1<k_2\cdots<k_n.$$
Generically, we have $k_j=j+1$ for every $j\in[\![1,n]\!]$. This is the reason why we consider from now on the following scalar conservation law, which we call the {\em multi-dimensional Burgers equation}~:
\begin{equation}
\label{eq:Bur}
\partial_tu+\partial_j\frac{u^2}2\,+\cdots+\,\partial_n\frac{u^{n+1}}{n+1}=0.
\end{equation}
This particular flux was already considered by G. Crippa et al. \cite{COW}.
If $n=1$, we recognize the original Burgers equation. The equation (\ref{eq:Bur}) is a prototype for genuinely nonlinear conservation laws, those which satisfy the assumption
\begin{equation}
\label{eq:GNL}
\det(f'',\ldots,f^{(n+1)})\ne0.
\end{equation}
The latter condition is a variant of the {\em non-degeneracy condition} at work in the kinetic formulation of the equation (\ref{eq:mdcl})~; see \cite{LPT} or \cite{Per}.

\bigskip

Let us review two preliminary answers to our natural question, in the context of (\ref{eq:Bur}).
\begin{itemize}
\item On the one hand, we might assume that $u_0\in L^1\cap L^p(\R^n)$ for some $p>1$. Let us define $$u_{0m}:=\max(-m,\min(u_0,m))\in L^1\cap L^\infty(\R^n),$$ 
which tends towards $u_0$ in the $L^1$-norm. We have $u=\lim_{m\rightarrow+\infty}u_m$, where $u_m$ is the solution associated with the data $u_{0m}$, and the limit stands in $C_b(\R_+;L^1(\R^n))$. Because of $\|u_m(t)\|_p\le\|u_{0m}\|_p\le\|u_0\|_p$, the sequence $(u_m)_{m\ge1}$ is bounded in $L^\infty(\R_+;L^p(\R^n))$. We infer that $u\in L^\infty(\R_+;L^p(\R^n))$, and $u_m$ converges towards $u$ in $L^\infty(\R_+;L^q(\R^n))$ for every $q\in[1,p)$. In addition, $u_m$ converges weakly in $L^\infty(\R_+;L^p(\R^n))$. If $p>n+1$, we may pass to the limit as $m\rightarrow+\infty$ in the sequences
$$u_m^k,\qquad {\rm sgn}(u_m-a)(u_m^k-a^k).$$
Passing to the limit in the integral formulations (\ref{eq:distr}) and (\ref{eq:entrina}), we conclude that $u$ is a genuine entropy solution of the Cauchy problem. Notice that the argument does not work out when $p=n+1$, because of the last component of the flux: we are not certain that $u_m^{n+1}$ converges in $L^1_{\rm loc}$ towards $u^{n+1}$. If $p>n+2$, we find as well that $u$ satisfies the entropy inequality for the pair $(\bar\eta,\bar q)$.

The drawback of this argument is that it does not exploit the nonlinearity of the equation, a property which is expected to imply some kind of regularization or dispersion (see Theorem 4 and Proposition 1 of \cite{LPT}). We should be able to lower somehow the threshold $p>n+1$.
\item The other answer concerns the one-dimensional case ($n=1$). The Kruzhkov solution of the classical Burgers equation satisfies the inequality
\begin{equation}
\label{eq:Daf}
TV\left(\frac{u(t)^2}2\right)\le\frac{2\|u_0\|_1}t\,,
\end{equation}
due to B\'enilan \& Crandall \cite{BC}, who exploit the homogeneity of the flux. It is extended by Dafermos \cite{Daf} to situations where the flux $f$ has an inflexion point and the data $u_0$ has bounded variations, by a careful use of the generalized backward characteristics.  It implies in particular an estimate
\begin{equation}
\label{eq:heat}
\|u(t)\|_\infty\le2\,\sqrt{\frac{2\|u_0\|_1}t\,}\,.
\end{equation}
This shows that the assumption $u_0\in L^1(\R)$ is sufficient in order that $u$ be a true entropy solution. This is definitely better than the threshold $L^1\cap L^2(\R^n)$ considered in the previous paragraph. 

Dafermos' argument, which is the most general one, uses the ordered structure of the real line. Backward characteristics are not unique in general. Given a base point $(x^*,t^*)$ in the upper half-plane, one has to define and analyse the minimal and the maximal ones. These notions have not yet been extended to the multi-dimensional situation (see however \cite{Sil} for a weaker notion).
\end{itemize}

Our main result here is the following statement. It tells us that $L^1(\R^n)$ is the right space for initial data.
\begin{thm}[Multi-d Burgers equation.]\label{th:wpd}
Suppose $u_0\in L^1(\R^n)$. Define $u(t)=S_tu_0$ and set  $u(t,y)=u(t)(y)$ for $t>0$ and $y\in\R^n$. Then 
\begin{enumerate}
\item There holds an algebraic decay:
\begin{equation}
\label{eq:decay}
\|u(t)\|_{\frac{d^2}{d-1}}\le c_dt^{-\delta}\,\|u_0\|_1^\gamma ,
\end{equation}
where
$$\gamma=\frac{d^2+1}{d(d^2-d+2)}\,,\qquad\delta=2\,\frac{(d-1)(d^2-d+1)}{d^2(d^2-d+2)}\,<1.$$
\item For every $k\in[\![1,d+1]\!]$, there holds $u^k\in L^1_{\rm loc}(\R_+\times\R^n)$.
\item The function $u$ is an entropy solution of the Cauchy problem.
\item It satisfies the additional entropy inequality
$$\partial_t\bar\eta(u)+{\rm div}_y\bar q(u)\le0.$$
\item If in addition $u_0\in L^1\cap L^d(\R^n)$, then $u\in L^{\frac{d^2}{d-1}}(\R_+\times\R^n)$ and there holds
\begin{equation}
\label{eq:estfond}
\left(\int_0^\infty\!\int_{\R^n}u^{\frac{d^2}{d-1}}dydt\right)^{\frac{d-1}d}\le c_d\left(\int_{\R^n}u_0(y)^ddy\right)^{\frac12}\left(\int_{\R^n}u_0(y)dy\right)^{\frac12}.
\end{equation}
\end{enumerate}
\end{thm}

\paragraph{Comments.}
\begin{itemize}
\item The assumption that $u_0\in L^1(\R^n)$ extends that available in the $1$-dimension situation. However, when $n=1$, Theorem \ref{th:wpd} provides an estimate of $u(t)$ in $L^4(\R)$ only, instead of the known $L^\infty(\R)$ or $BV(\R)$. Our results are new only when $n\ge2$.
\item The decay result is optimal when $n=1$, where it states that
$$\|u(t)\|_4\le {\rm cst}\,\cdot t^{-3/8}\|u_0\|_1^{5/8}.$$
This is the exact rate for an N-wave
$$N_L(t,y)=\left\{\begin{array}{lcr}
\frac{y}{1+t} & \hbox{if} & y\in(0,L\sqrt{1+ t}), \\
0 & \hbox{otherwise.} &
\end{array}\right.$$
It raises therefore the question whether the decay rate given by (\ref{eq:decay}) is accurate also when $n\ge2$. 
\item Estimate (\ref{eq:estfond}) ressembles a Strichartz inequality. It seems to be new in this situation where the principal part in not a linear operator, but a quasilinear one.
\item By H\"older interpolation, together with $u\in L^\infty(\R_+;L^1(\R^n))$, (\ref{eq:decay}) implies
\begin{equation}\label{eq:gendec}
\|u(t)\|_q\le c_{d,q}t^{-\kappa/q'}\|u_0\|_1^{1-\nu/q'},\qquad\forall q\in(1,\frac{d^2}{d+1}),
\end{equation}
where 
$$\kappa=2\,\frac{d-1}{d^2-d+2}\,\quad\hbox{and}\quad\nu=\frac{d(d-1)}{d^2-d+2}\,.$$
\item A natural question is whether (\ref{eq:gendec}) extends to every $q\in[1,\infty]$. In particular, is it true that $u(t)\in L^\infty(\R^n)$ for every $t>0$~? We now that it is true when $n=1$, see (\ref{eq:heat}).
\item
A useful contribution in this direction was obtained recently by L. Silvestre \cite{Sil}, whose Theorem 1.5 tells in particular that if $u_0\in L^1\cap L^\infty(\R^n)$, then
$$\|u(t)\|_\infty\le C(\|u_0\|_\infty;\mu)\,\|u_0\|_1^\mu\,t^{-n\mu}$$
for every $\mu<\mu_0$ where
$$\mu_0=\frac2{d^2-d+2}\,.$$
This decay is almost the same as that suggested by extrapolation to $q=\infty$ of ours, because of
$$n\mu_0=\kappa=\frac\kappa{q'}\,.$$
It would be exactly that one if the limit exponent $\mu_0$ was allowed, and the dependency of the constant upon $\|u_0\|_\infty$ was removed.
\end{itemize}

\subsection{Other ``monomial'' scalar conservation laws}

As suggested above, we may be interested into more general conservation laws, whose fluxes are monomial. Denoting $P_k(s)=\frac{s^k}k$\,, consider the PDE
\begin{equation}
\label{eq:monom}
\partial_tu+\partial_1P_{k_1}(u)+\cdots+\partial_nP_{k_n}(u)=0,
\end{equation}
where $1<k_1<\cdots<k_n$ are integers. We leave, as a tedious exercise, the interesting reader to adapt the calculations of the two next sections to (\ref{eq:monom}), to prove the following result. We denote
$$K=\sum_{j=1}^nk_j,\qquad N=1+2K-n.$$
\begin{thm}\label{th:mon}
Suppose that $nk_n<N$. Then for every initial data $u_0\in L^1(\R^n)$, the abstract solution given by the continuous extension of the semi-group $(S_t)_{t\ge0}$ to $L^1(\R^n)$, is actually an entropy solution of the Cauchy problem for (\ref{eq:monom}). It satisfies a dispersion estimate
$$\left(\int_0^\infty dt\int_{\R^n}u(t,y)^{\frac Nn}dy\right)^{\frac n{n+1}}\le c_d\left(\int_{\R^n}u_0(y)dy\right)^{1-\theta}\left(\int_{\R^n}u_0(y)^{k_n}dy\right)^{\theta},$$
where 
$$\theta=\frac{K-n}{(n+1)(k_n-1)}\,.$$
It decays as follows
$$\|u(t)\|_{\frac Nn}\le c_d\frac{\|u_0\|_1^\gamma}{t^\delta}$$
where
$$\gamma= \frac nN+\frac{N-n}{N(N-K)}\,,\qquad \delta=\frac{n(N-n)}{N(N-K)}\,.$$
\end{thm}

\bigskip

The r\^ole of the assumption $nk_n<N$ is to allow us to estimate $\|u(t)\|_{k_n}$ in terms of $\|u(t)\|_1$ and $\|u(t)\|_{N/n}$, in order to apply a Gronwall argument to the dispersion estimate. Notice that it is always satisfied in one space dimension, because then $1\cdot k_1<N=2k_1$

\paragraph{Remark.} If $k_n$ is larger than $\frac Nn$\,, there should be a weaker result. There will be some exponent $p=p(k_n,N)\in(1,k_n)$ such that if $u_0\in L^1\cap L^p(\R^n)$, then the abstract solution is actually an entropy solution. We leave the calculation of $p(k_n,N)$ to the motivated reader.

\section{Proof of Estimate (\protect\ref{eq:estfond})}\label{s:disp}

Because $u$ is obtained as the limit in $C(\R_+;L^1(\R^n))$ of $u_m$, the solution associated with the data $u_{0m}={\rm Proj}_{[-m,m]}u_0$, the estimates (\ref{eq:decay}) and (\ref{eq:estfond}) need only to be proved when the initial data belongs to $L^1\cap L^\infty(\R^n)$, that is within Kruzhkov's theory. Then they extend to $L^1$-data by a density argument. 

When $u_0\in L^1\cap L^\infty$, (\ref{eq:estfond}) will provide a uniform bound of
$$\int_0^\infty\|u_m(t)\|_{\frac{d^2}{d-1}}^{\frac{d^2}{d-1}}dt.$$
Then, because of $u_m\rightarrow u$ in $C(\R_+;L^1(\R^n))$ as $m\rightarrow+\infty$, we infer by interpolation that the convergence holds true in every space $L^q(\R_+;L^p(\R^n))$ for which
$$1\le p<\frac{d^2}{d-1}\quad(\hbox{equivalently }\frac{d^2}{d-1}<q\le\infty)\qquad\hbox{and}\quad\frac{d-1}p+\frac{d^2-d+1}q=d-1.$$

\bigskip

Because of $(S_tu_0)_\pm\le S_t(u_0)_\pm$,  it is enough to consider data that are either non-negative or non-positive. But since $v(t,y)=-u(t,-y_1,y_2,\ldots,(-1)^ny_n)$ is the entropy solution associated with $v_0(y)=-u_0(-y_1,y_2,\ldots,(-1)^ny_n)$, it suffices to prove (\ref{eq:estfond}) for non-negative data and solutions. We therefore assume from now on that $u_0\ge0$, and thus $u\ge0$ over $\R_+\times\R^n$.

\bigskip

If $a\in\R$, we define a symmetric matrix
$$M(a)=\left(\frac{a^{i+j+1}}{i+j+1}\right)_{0\le i,j,\le n}.$$
Remarking that
$$M(a)=\int_0^aV(s)\otimes V(s)\,ds,\qquad V(s)=\begin{pmatrix} 1 \\ \vdots \\ s^n \end{pmatrix},$$
we obtain that $M(a)$ is positive definite whenever $a>0$. We have obviously
$$\det M(a)=H_d\,a^{d^2},$$
where 
$$H_d=\left\|\frac1{i+j+1}\right\|_{0\le i,j,\le n}>0$$
is the determinant of the Hilbert matrix (this is the only case where we do not write $c_d$ for a dimensional constant).

Let us form the symmetric tensor
$$T(t,y)=M(u(t,y)),$$
with positive semi-definite values. Its first line is formed of $(u,f(u))$ and therefore is divergence-free by (\ref{eq:Bur}). The second line is formed of $(\bar\eta(u),\bar q(u))$, an entropy-flux pair. It is not divergence-free in general, although it is so away from shock waves and other singularities of the solution $u$. But the entropy inequality tells us that the opposite of its divergence if a non-negative, hence bounded measure, 
$$\mu_1=-{\rm div}_{t,y}(\bar\eta(u),\bar q(u))\ge0.$$
The total mass of $\mu_1$ over a slab $(0,\tau)\times\R^n$ is given by
$$\|\mu_1\|=\int_{\R^n}\bar\eta(u_0(y))\,dy-\int_{\R^n}\bar\eta(u(\tau,y))\,dy\le\int_{\R^n}\bar\eta(u_0(y))\,dy.$$
Notice that the latter bound does not depend of $\tau$.
The same situation occurs for the other lines of $T$. They are of the form $(\eta(u),q(u))$ where $(\eta,q)$ is an entropy-flux pair with $\eta$ convex over $\R_+$ (recall that $u$ takes only non-negative values). The distribution
$$\mu=-{\rm div}_{t,y}(\eta(u),q(u))$$ is therefore again a bounded measure, whose total mass over $\R_+\times\R^n$ is bounded by 
$$\int_{\R^n}\eta(u_0(y))\,dy.$$

We conclude that the row-wise divergence of $T$ is a (vector-valued) bounded measure, whose total mass is bounded above by
$$\sum_{j=2}^d\int_{\R^n}\frac{u_0(y)^j}j\,dy.$$
We may therefore apply Compensated integrability (Theorems 2.2 and 2.3 of \cite{Ser_CI}) to the tensor $T$, that is
$$\int_0^\tau dt\int_{\R^n}(\det T)^{\frac1{d-1}}dy\le c_d\left(\|T_{0\bullet}(0,\cdot)\|_1+\|T_{0\bullet}(\tau,\cdot)\|_1+\|{\rm Div}_{t,y}T\|_{{\cal M}((0,\tau)\times\R^n)}\right)^{\frac d{d-1}}.$$
Because of
$$\|T_{0\bullet}(t,\cdot)\|_1\le\sum_{j=1}^d\int_{\R^n}\frac{u(t,y)^j}j\,dy\le\sum_{j=1}^d\int_{\R^n}\frac{u_0(y)^j}j\,dy,$$
we deduce
\begin{equation}
\label{eq:nonhom}
\int_0^\tau dt\int_{\R^n}u^{\frac{d^2}{d-1}}dy\le c_d\left(\sum_{j=1}^d\int_{\R^n}u_0(y)^j\,dy\right)^{\frac d{d-1}}.
\end{equation}

\bigskip

The only bad feature in the estimate (\ref{eq:nonhom}) is the lack of homogeneity of its right-hand side. To recover a well-balanced inequality, we exploit an idea already used in \cite{Ser_DPT}. We begin by remarking that if $\lambda>0$ is a constant parameter, then the function
$$v(t,y)=\frac1\lambda\,u(\lambda t,\lambda^2y_1,\ldots,\lambda^dy_n)$$
is the entropy solution associated with the initial data
$$v_0(y)=\frac1\lambda\,u_0(\lambda^2y_1,\ldots,\lambda^dy_n).$$
Applying (\ref{eq:nonhom}) to the pair $(v,v_0)$ instead, then using
$$\int_{\R^n}v_0(y)^jdy=\lambda^{-j-\frac{(d-1)(d+2)}2}\int_{\R^n}u_0(y)^jdy$$
and 
$$\left(\int_0^\infty\!\int_{\R^n}v^{\frac{d^2}{d-1}}dydt\right)^{\frac{d-1}d}=\lambda^{-d-\frac{d^2-1}2}\left(\int_0^\infty\!\int_{\R^n}u^{\frac{d^2}{d-1}}dydt\right)^{\frac{d-1}d},$$ 
we get a parametrized inequality
$$\left(\int_0^\tau dt\int_{\R^n}u^{\frac{d^2}{d-1}}dy\right)^{\frac{d-1}d}\le c_d\lambda^{\frac{d+1}2}\sum_{j=1}^d\lambda^{-j}\int_{\R^n}u_0(y)^j\,dy.$$
In order to minimize the right-hand side, we choose the  value
$$\lambda=\left(\int_{\R^n}u_0(y)^ddy/\int_{\R^n}u_0(y)dy\right)^{\frac1{d-1}}.$$
The extreme terms, for $j=1$ or $d$, contribute on a equal foot with
$$\left(\int_{\R^n}u_0(y)^ddy\right)^{\frac12}\left(\int_{\R^n}u_0(y)^2dy\right)^{\frac12}.$$
The other ones, which are
$$\left(\int_{\R^n}u_0(y)^ddy/\int_{\R^n}u_0(y)dy\right)^{\frac{d+1-2j}{2(d-1)}}\int_{\R^n}u_0^jdy,$$
are bounded by the same quantity, because of H\"older inequality. We end therefore with the fundamental estimate (\ref{eq:estfond}).

\section{Proof of Theorem \protect\ref{th:wpd}}\label{s:dec}

We now complete the proof of our main theorem.

\subsection{The decay result}

We keep working with the assumptions $u_0\in L^1\cap L^\infty(\R^n)$ and $u_0\ge0$.

\bigskip

Let us  define
$$X(t):=\int_{\R^n}u^{\frac{d^2}{d-1}}dy=\|u(t)\|_{\frac{d^2}{d-1}}^{\frac{d^2}{d-1}}.$$
From the H\"older inequality, we have
$$\|u_0\|_d\le \|u_0\|_1^{1-\theta}\|u_0\|_{\frac{d^2}{d-1}}^\theta,\qquad\theta=\frac{d(d-1)}{d^2-d+1}\,\in(0,1).$$
The inequality (\ref{eq:estfond}) implies therefore
$$\int_0^\infty X(t)\,dt\le c_d \|u_0\|_1^\alpha X(0)^\beta$$
where
$$\beta=\frac\theta2\,,\qquad\alpha=\frac{d(d^2+1)}{2(d-1)(d^2-d+1)}\,.$$
Considering the solution $v(t,y)=u(t+\tau,y)$, whose initial data is $u(\tau,\cdot)$, we also have
\begin{equation}
\label{eq:Xtau}
\int_\tau^\infty X(t)\,dt\le c_d \|u(\tau)\|_1^\alpha X(\tau)^\beta\le c_d \|u_0\|_1^\alpha X(\tau)^\beta.\end{equation}
Let us denote 
$$Y(\tau):=\int_\tau^\infty X(t)\,dt.$$
We recast (\ref{eq:Xtau}) as
$$Y^{\frac1\beta}+c_d\|u_0\|_1^{\frac\alpha\beta}Y'\le0.$$
Multiplying by $Y^{-1/\beta}$ and integrating, we infer (mind that $1-\frac1\beta$ is negative)
$$t+c_d\|u_0\|_1^{\frac\alpha\beta}Y(0)^{1-\frac1\beta}\le c_d\|u_0\|_1^{\frac\alpha\beta}Y(t)^{1-\frac1\beta}.$$
This provides a first decay estimate
$$Y(t)\le c_d\|u_0\|_1^{\frac\alpha{1-\beta}}t^{-\frac\beta{1-\beta}}.$$
Remarking that $t\mapsto X(t)$ is a non-increasing function, so that
$$\frac\tau2\,X(\tau)\le Y(\frac\tau2),$$
we deduce the ultimate decay result
$$X(t)\le c_d\|u_0\|_1^{\frac\alpha{1-\beta}}t^{-\frac1{1-\beta}}.$$
Restated in terms of a Lebesgue norm of $u(t)$, it says
\begin{equation}
\label{eq:decayi}
\|u(t)\|_{\frac{d^2}{d-1}}\le c_d\|u_0\|_1^\gamma\, t^{-\delta},
\end{equation}
where
$$\gamma=\frac{d^2+1}{d(d^2-d+2)},\qquad\delta=2\,\frac{(d-1)(d^2-d+1)}{d^2(d^2-d+2)}\,.$$

\subsection{The function $u$ is an entropy solution}

We already know that the functions $u_m$ are entropy solutions, with initial data $u_{0m}\in L^1\cap L^\infty(\R^n)$. Because of (\ref{eq:estfond}), we have seen that $u_m$ converges towards $u$ in the norm of $L^q(\R_+;L^p(\R^n))$ whenever
$$1\le p<\frac{d^2}{d-1}\quad\hbox{and}\quad\frac{d-1}p+\frac{d^2-d+1}q=d-1.$$
Because of $\frac{d^2}{d-1}>d+1$, this implies that $(u_m)^k\rightarrow u^k$ for every $k\in[\![1,d+1]\!]$. Recalling that
$$f(s)=(s,\frac{s^2}2\,,\ldots,\frac{s^d}d)\quad\hbox{and}\quad \bar q(s)=(\frac{s^3}3\,,\ldots,\frac{s^{d+1}}{d+1}),$$
we see that we may pass to the limit as $m\rightarrow+\infty$ in the identity
$$\int_0^\infty dt\int_{\R^n}(u_m\partial_t\phi+f(u_m)\cdot\nabla_y\phi)\,dy+\int_{\R^n}u_{0m}(y)\phi(0,y)\,dy=0,\qquad\forall\phi\in{\cal D}(\R^{1+n}),$$
as well as in the Kruzhkov inequalities
$$\int_0^\infty dt\int_{\R^n}(\eta_a(u_m)\partial_t\phi+q_a(u_m)\cdot\nabla_y\phi)\,dy+\int_{\R^n}\eta_a(u_{0m}(y))\phi(0,y)\,dy\ge0,\qquad\forall\phi\in{\cal D}^+(\R^{1+n}),$$
and in the inequality
$$\int_0^\infty dt\int_{\R^n}(\bar\eta(u_m)\partial_t\phi+\bar q(u_m)\cdot\nabla_y\phi)\,dy+\int_{\R^n}\bar\eta(u_{0m}(y))\phi(0,y)\,dy\ge0,\qquad\forall\phi\in{\cal D}^+(\R^{1+n}).$$
Therefore $u$ is an entropy solution with initial data $u_0$, which satisfies in addition the entropy inequality for the pair $(\bar\eta,\bar q)$.

\paragraph{Remark.} When $n\ge2$, the Compensated Integrability cannot be applied directly to the solution $u$, when the data is only integrable. Because we don't know whether the $j$th line of $T$ is locally integrable if $j=3,\cdots,n+1$~; its last component is $\frac{u^{n+j}}{n+j}$\,, where the exponent $n+j$ is larger than $\frac{d^2}{d-1}$\,.

\section{The strategy for general fluxes $f$}\label{s:gen}

We come back to the study of a multi-dimensional conservation law of general form (\ref{eq:mdcl}). Following the ideas develloped in the Burgers case, we begin by considering a signed, bounded initial data: $u_0\in L^1\cap L^\infty(\R^n)$, $u_0\ge0$. If $a\in\R_+$, we define a symmetric matrix
$$M(a)=\int_0^aF'(s)\otimes F'(s)\,ds,$$
where $F(s)=(s,f_1(s),\ldots,f_n(s))$. This matrix is positive definite under the non-degeneracy condition that $F([0,a])$ is not contained in an affine hyperplane. We denote
$$\Delta(a):=(\det M(a))^{\frac1n}\ge0.$$

Let us define $T(t,y):=M(u(t,y))$. Because $u\in L^\infty(0,\tau;L^1\cap L^\infty(\R^n))$, the tensor $T$ is integrable over $(0,\tau)\times\R^n$.
The first line of $T$ is divergence-free. The other lines are made of entropy-entropy flux pairs $(f_i,Q_i)$. Since $f_i$ might not be convex, we cannot estimate the measure $\mu_i=-\partial_tf_i(u)-{\rm div}_yQ_i(u)$ directly by the integral of $f_i(u_0)$. To overcome this difficulty, we define a convex function $\phi$ over $\R_+$ by
$$\phi(0)=\phi'(0)=0,\qquad\phi''(s)=|f''(s)|.$$
Remark that $|f'|\le\phi'$ and $|f|\le\phi$. Let $\Phi$ be the entropy flux associated with the entropy $\phi$. Then the measure $\nu:=-\partial_t\phi(u)-{\rm div}_y\Phi(u)$ is non-negative and a bound of its total mass is  as usual
$$\|\nu\|\le\int_{\R^n}\phi(u_0(y))\,dy.$$

We now use the kinetic formulation of (\ref{eq:mdcl}), a notion for which we refer to \cite{Per}, Theorem 3.2.1. Recall the definition of the kinetic function $\chi(\xi;a)$, whose value is ${\rm sgn}\,a$ if $\xi$ lies between $0$ and $a$, and is $0$ otherwise. There exists a non-negative bounded measure $m(t,y,\xi)$ such that the function $g(t,y,\xi)=\chi(\xi;u(t,y))$ satisfies
$$\partial_tg+f'(\xi)\cdot\nabla_yg=\frac{\partial}{\partial\xi}\,m,\qquad g(0,y;\xi)=\chi(\xi;u_0(y)).$$
If $(\eta,q)$ is an entropy-entropy flux pair, then the measure $\mu=-\partial_t\eta-{\rm div}_yq$ is given by
$$\mu=\int_\R \eta''(\xi)dm(\xi).$$
We deduce that the vector-valued measure $\mu=(\mu_1,\ldots,\mu_n)$ satisfies $|\mu|\le\nu$. This yields the estimate
$$\|\mu\|\le\int_{\R^n}\phi(u_0(y))\,dy.$$
We may therefore apply the compensated integrability, which gives here
$$\int_0^\tau dt\int_{\R^n}\Delta(u(t,y))\,dy\le c_d\left(\|F(u_0)\|_1+\|F(u(\tau))\|_1+\int_{\R^n}\phi(u_0(y))\,dy\right)^{1+\frac1n}.$$
Because of $|f|\le\phi$ and $\|\phi(u(\tau))\|_1\le\|\phi(u_0)\|_1$, we end up with an analog of (\ref{eq:estfond})
\begin{equation}
\label{eq:CIgen}
\int_0^\infty dt\int_{\R^n}\Delta(u(t,y))dy\le c_d(\|u_0\|_1+\|\phi(u_0)\|_1)^{1+\frac1n}.
\end{equation}

\bigskip

To improve the inequality above, we use again a scaling argument. However, because the components $f_j$ of the flux are not homogeneous anymore, we modify simultaneously the solution and the flux, using the fact that the constant $c_d$ in (\ref{eq:CIgen}) does not depend upon $f$. Our new dependent variables are
$$v(t,y)=\frac1\lambda\,u(\lambda t,Py),\qquad v_0(y)=v(0,y)=\frac1\lambda\,u_0(Py)$$
where $P\in {\bf GL}_n(\R)$ is a matrix to be chosen later. The function $v$ is an entropy solution of the Cauchy problem associated with  the conservation law
$$\partial_tv+{\rm div}_y\,g(v)=0$$
for the flux
\begin{equation}
\label{eq:clmod}
g(s)=P^{-1}f(\lambda s).
\end{equation}
The symmetric matrix $N(a)$ that plays the role of $M(a)$ for (\ref{eq:clmod}) is given by the formula 
$$N(s)=\lambda\, RM(\lambda s)R^T,\qquad R={\rm diag}(\frac1\lambda,P^{-1}).$$
We have $\det N(s)=\lambda^{n-1}(\det P)^{-2}\det M(\lambda s)$, from which we derive
$$\int_0^\infty\!\int_{\R^n}(\det N(v))^{\frac1n}dy\,dt=\lambda^{-\frac1n}(\det P)^{-1-\frac2n}\int_0^\infty\!\int_{\R^n}\Delta(u)\,dy\,dt.$$
When applying (\ref{eq:clmod}) to $v$ and $g$, one integral in the right-hand side transforms easily:
$$\int_{\R^n}v_0(y)\,dy=\frac1{\lambda\det P}\,\int_{\R^n}u_0(y)\,dy.$$
The other one involves  a modified function
$$\phi_{\lambda,P}(s)=\psi_P(\lambda s),\qquad\psi_P''(s)=|P^{-1}f''(s)|.$$
As usual $\psi_P$ is fixed by $\psi_P(0)=\psi_P'(0)=0$. We have therefore
$$\int_{\R^n}\phi_P(v_0(y))\,dy=\frac1{\det P}\,\int_{\R^n}\psi_P(u_0(y))\,dy.$$
All the identities above, together with (\ref{eq:clmod}) applied to $(v,g)$, yield our parametrized estimate
$$\int_0^\infty\!\int_{\R^n}\Delta(u)\,dy\,dt\le c_d\lambda^{\frac1n}(\det P)^{\frac1n}\left(\frac1\lambda\,\int_{\R^n}u_0dy+\int_{\R^n}\psi_P(u_0)\,dy\right)^{1+\frac1n}.$$
We optimize this inequality with respect to $\lambda$, by choosing
$$\lambda=\frac{\int_{\R^n}u_0dy}{\int_{\R^n}\psi_P(u_0)\,dy}.$$
We infer
$$\int_0^\infty\!\int_{\R^n}\Delta(u)\,dy\,dt\le c_d(\det P)^{\frac1n}\left(\int_{\R^n}u_0dy\right)^{\frac1n}\int_{\R^n}\psi_P(u_0)\,dy.$$
There remains to minimize the right-hand side with respect to $P$~:
\begin{equation}
\label{eq:Grongen}
\int_0^\infty\!\int_{\R^n}\Delta(u)\,dy\,dt\le c_d\left(\int_{\R^n}u_0dy\right)^{\frac1n}I[u_0],
\end{equation}
where
$$I[w]:=\inf_{P\in{\bf GL}_n(\R)}(\det P)^{\frac1n}\int_{\R^n}\psi_P(w(y))\,dy.$$
The calculation of $I[w]$ has to be made on a case-by-case basis. 

\bigskip

Let us define again
$$X(t)=\int_{\R^n}\Delta(u(t,y))\,dy.$$
Applying (\ref{eq:Grongen}) on an interval $(\tau,\infty)$ instead, and using the decay of the $L^1$-norm, we arrive to
$$\int_\tau^\infty X(t)\,dt\le c_d\left(\int_{\R^n}u_0dy\right)^{\frac1n}I[u(\tau)].$$
A decay result will be obtained through a Gronwall argument, whenever we can estimate $I[u(t)]$ in terms of $\|u(t)\|_1$ and $X(t)$.

\end{document}